\documentclass{article}
\usepackage[utf8]{inputenc}
\usepackage{amsmath,amssymb}
\usepackage{algorithmic}
\usepackage[ruled,commentsnumbered,linesnumbered,vlined]{algorithm2e}
\usepackage{xcolor}
\usepackage{indentfirst}
\usepackage{comment}
\usepackage{longtable}
\usepackage{lscape}
\usepackage{rotating}

\usepackage{pdflscape}

\DeclareMathOperator*{\argmin}{argmin}

\usepackage{lscape}
\usepackage{longtable}
\usepackage{amsthm}
\usepackage{subfigure}
\usepackage{picture}
\usepackage{graphicx}
\usepackage{fullpage}
\usepackage{hyperref}
\usepackage{float}

\usepackage{apacite}

\DeclareMathAlphabet\mathbfcal{OMS}{cmsy}{b}{n}

\newcommand{\rafaelC}[1]{\textcolor{black}{#1}}

\title{A fast and effective MIP-based heuristic for a selective and periodic inventory routing problem in reverse logistics}
\author{ 
Leopoldo E. Cárdenas-Barrón {\thanks{Tecnologico de Monterrey, School of Engineering and Sciences,  Ave. Eugenio Garza Sada 2501, Monterrey, N.L., Mexico, 64849. ({\tt lecarden@tec.mx})}}
    \and
    Rafael A. Melo {\thanks{Universidade Federal da Bahia, Departamento de Ci\^{e}ncia da Computa\c{c}\~{a}o, Computational Intelligence and Optimization Research Lab (CInO), Salvador, Brazil.  ({\tt melo@dcc.ufba.br}) }}
       }

\begin{document}

\maketitle

\begin{abstract}

    We consider an NP-hard selective and periodic inventory routing problem (SPIRP) in a waste vegetable oil collection environment. 
    This SPIRP arises in the context of reverse logistics where a biodiesel company has daily requirements of oil to be used as raw material in its production process. These requirements can be fulfilled by using the available inventory, collecting waste vegetable oil or purchasing virgin oil. The problem consists in determining a period (cyclic) planning for the collection and purchasing of oil such that the total collection, inventory and purchasing costs are minimized, while meeting the company's oil requirements and all the operational constraints.
 We propose a MIP-based heuristic which solves a relaxed model without routing, constructs routes taking into account the relaxation’s solution and then improves these routes by solving the capacitated vehicle routing problem associated to each period. Following this approach, an \textit{a posteriori} performance guarantee is ensured, as the approach provides both a lower bound and a feasible solution.  The performed computational experiments show that the MIP-based heuristic is very fast and effective as it is able to encounter near optimal solutions with low gaps within seconds, improving several of the best known results using just a fraction of the time spent by a state-of-the-art heuristic. A remarkable fact is that the proposed MIP-based heuristic improves over the best known results for all the large instances available in the literature. \\

    \noindent \textbf{Keywords:} reverse logistics; inventory routing; mixed integer programming; heuristics; sustainability.
\end{abstract}

\section{Introduction}
\label{sec:introduction}

Nowadays, several industries have an increasing concern about the implementation of green actions in order to collaborate with the sustainability of the environment. For this reason, \citeA{LinChoHoChuLam14} argued that green logistics is a big challenge that all companies must consider in their agendas. An important area of study that deals with sustainable supply chains is reverse logistics~\cite{DekFleIndWas13}, as it includes the application of methods which take into consideration, among other characteristics, the reuse of recyclable material.

The inventory routing problem (IRP) is a variant of the vehicle routing problem (VRP) which integrates inventory and vehicle routing decisions in a sole formulation. It is important to mention that the IRP models numerous distribution and collection situations that occur in practice. \citeA{AndHenHofArilChr10} and \citeA{CoeCorLap14} presented good surveys with the aim of having a good comprehension of the IRP problem. On the one hand, \citeA{AndHenHofArilChr10} dealt with the IRP’s applications. On the other hand, \citeA{CoeCorLap14} studied IRP methodological aspects. \citeA{CoeLap13} solved the following three variants of the IRP by applying a branch-and-cut algorithm: the multi-vehicle IRP considering both homogeneous and heterogeneous fleet, the IRP with transshipment options, and the IRP with added consistency features. 
In a subsequent paper, \citeA{CoeLap14} built an exact mathematical formulation by implementing valid inequalities for some types of IRPs. \citeA{MjiJarMacHanMla14} considered a multi-product IRP and tackled it using a variable neighborhood search. \citeA{SoyBloHaiVor15} formulated a multi-period IRP for perishable products considering that demand is uncertain and \citeA{SoyBloHaiVor18} introduced a green IRP for perishable items. \citeA{CorLagMusVoc15} considered the multi-product IRP and developed a three-stage heuristic algorithm taking into account a decomposition of the seller’s decision process. \citeA{Raa15} proposed a family of cyclic IRP. The aforementioned papers are related to forward logistics. However, some authors also considered IRPs in reverse logistics. In this direction, \citeA{AksKaySalAkc12} formulated an IRP for a reverse logistics situation where a biodiesel company requires oil to be used as raw material, which can be obtained as both virgin oil and waste vegetable oil, in its production process. Thus, the collection of waste vegetable oil is a significant reverse operation to be performed with the aims of having economic gains and improving environmental sustainability. 

In fact, \citeA{AksKaySalAkc12} proposed a selective and periodic inventory routing problem (SPIRP) for a waste vegetable oil collection, which is the problem we consider in our work. In a subsequent paper, \citeA{AksKaySibOzg14} presented an adaptive large neighborhood search algorithm (ALNS) to solve the SPIRP. Recently, \citeA{CarGonTreGar19} developed a reduce and optimize approach (ROA) heuristic for this SPIRP. The ROA defines a reduced feasible region based on the original problem to be optimized. 
The ROA heuristic was able to outperform the ALNS of \citeA{AksKaySibOzg14} for most of the benchmark instances of the SPIRP, at the expense of computational time.

\rafaelC{Other variants of routing problems were also considered by several authors. \citeA{bertazzi2015} addressed a stochastic inventory routing problem with one supplier, multiple retailers, and transportation procurement in which the supplier has a constrained production capacity and the retailers’ demands are stochastic. The authors proposed a stochastic dynamic programming formulation that can be used to solve small instances to optimality. For large instances, they developed a matheuristic solution procedure that combines a rollout algorithm and mixed-integer linear programming models. \citeA{qiu2019} formulated a production inventory routing model with perishable inventory and proposed an exact branch and cut algorithm, which included some logical, strengthened lot-sizing and lifted Miller-Tucker-Zemlin subtour elimination inequalities. 
\citeA{neves2019} tackled a large multi-product production routing problem with delivery time windows using a three-phase procedure: firstly, the size of the original problem is reduced; secondly, the approach constructs an initial solution by decomposition; and thirdly, the initial solution is improved by solving different mixed-integer programming formulations. \citeA{fokkema2020} studied an inventory routing problem considering that decisions related to routing and inventory are made according to supply and not by demand. The authors introduced valid inequalities that allow to solve some instances to optimality. \citeA{alvarez2020} introduced a stochastic inventory routing problem for the case when both the supply and demand are stochastic. They proposed a solution procedure which is based on a progressive hedging algorithm.}

The main contribution of this paper is the proposal of a fast and effective MIP-based heuristic for the SPIRP. The proposed approach can provide near-optimal solutions whose optimality gaps are close to those obtained with a state-of-the-art approach, which consumes high computational time. As a matter of fact, computational experiments have proven the effectiveness of our approach, as new best results were encountered for several of the considered instances. To the best of our knowledge, this was a main open research area for the considered SPIRP.

The remainder of this paper is structured as follows. Section~\ref{sec:formulation} formally defines and describes a mathematical formulation of the SPIRP. Section~\ref{sec:mipheuristic} presents a new MIP-based heuristic for the SPIRP which basically solves a relaxation without routing variables and constraints, builds routes based on the relaxation’s solution, and then improves the routes using an approach for the capacitated vehicle routing problem (CVRP). Section~\ref{sec:experiments} presents the experiments carried out to evaluate the performance of the proposed MIP-based heuristic. Finally, section~\ref{sec:concludingremarks} gives some concluding remarks and future research avenues.










\section{Problem definition and mixed integer programming formulation}
\label{sec:formulation}

The selective and periodic inventory routing problem (SPIRP) can be formally defined as follows. 
There exist a set $IC =  \{1,\ldots,n\}$ of waste oil collection points (collection nodes) and one biodiesel production facility (depot node). The collection network can be seen as a complete directed graph in which there is an arc $(i,j)$ between every pair of nodes with an associated distance $d_{ij}$. A fixed cyclic planning horizon is defined with a predetermined set of periods $T=\{1,\ldots,\tau \}$ in each cycle. At each collection node $i \in IC$, a waste oil accumulation rate $a_{it}$ (in liters) at period $t\in T$ occurs. \rafaelC{To simplify the notation, it is assumed that the depot node has null waste oil accumulation rate for every period, i.e., $a_{0t}=0$ for $t\in T$.} \rafaelC{Additionally, define $A_i = \sum_{t \in T}a_{it}$ to be the total accumulation of waste oil at node $i\in IC$ during the planning horizon, which gives an upper bound on the amount of oil that should be collected from the node.} If a collection vehicle visits a collection node $i \in IC$ at period $t\in T$, then the whole waste oil stored up in period $t$ must be collected. This means that partial collection is not permitted. The biodiesel company requires $r_{t}$ liters of oil for each period $t\in T$ of the planning horizon and needs to determine, for each period, the amount of virgin oil to be purchased as well as the set of collection routes taking into consideration its oil requirements and the constraint that only one vehicle must be allocated to each route. Additional characteristics of the SPIRP are as follows. There is a fixed capacity $Q$ for each collection vehicle. The oil necessities of the biodiesel company can be covered from collected waste oil, purchased virgin oil, inventory on hand, or with a combination of these three sources. The costs involved in the SPIRP are: traveling cost $c$ per unit distance traveled, operating cost $v$ per vehicle per period, holding cost $h$ per liter of waste oil per period, and purchase cost $p$ per liter of virgin vegetable oil. The main goal is to decide how much virgin vegetable oil to buy and how much waste vegetable oil to gather from a set of collection points. Moreover, it is necessary to take decisions related to which collection nodes to include in the collection program, how many vehicles are needed, and their specific routes to define a periodic collection schedule that must be repeated in each cycle in such a manner that the total cost is minimized subject to the production requirements and vehicle capacity constraints.

\subsection{Mixed integer programming formulation}

In this section, we describe the formulation proposed in \citeA{AksKaySalAkc12}. Let the \textbf{index sets} be:
\begin{itemize}
    \item[--] $I = \{0,1,\ldots,n\}$;
    \item[--] $IC =  \{1,\ldots,n\}$;
    \item[--] $T=\{1,\ldots,\tau \}$.
\end{itemize}
Consider the \textbf{parameters}:
\begin{itemize}
    \item[--] $c$: traveling cost per unit distance;
    \item[--]  $d_{ij}$: distance from node $i$ to node $j$;
    \item[--]  $a_{it}$: amount of waste oil accumulated at node $i$ in period $t$;
    \item[--]  $r_t$: oil requirement in period $t$;
    \item[--]  $h$: per liter inventory holding cost at the depot;
    \item[--]  $v$: per vehicle fixed operating cost;
    \item[--]  $p$: per liter purchase cost of virgin vegetable oil;
    \item[--]  $Q$: vehicle capacity in liters\rafaelC{;}
    \item[--] \rafaelC{$A_i$: total amount of waste oil accumulated at node $i$ during the planning horizon, i.e., $A_i=\sum_{t \in T} a_{it}$.}
\end{itemize}
Define the \textbf{decision variables}:
\begin{itemize}
    \item[--]  $X_{ijt}$: binary variable which indicates whether arc $(i,j)$ is traversed by a vehicle in period $t$;
    \item[--]  $Y_{it}$: binary variable which indicates whether node $i$ is visited in period $t$;
    \item[--]  $Z_i$: binary variable which indicates whether node $i$ is visited at all during the planning cycle;
    \item[--]  $F_{ijt}$ the amount of waste oil traversing arc $(i,j)$ in period $t$;
    \item[--]  $W_{it}$: the amount of waste oil collected from node $i$ in period $t$;
    \item[--]  $I_{it}$: the amount of wast oil held as stock in node $i$ at the end of period $t$;
    \item[--]  $I_{i0}$: amount of waste oil in storage at the beginning of the planning cycle;
    \item[--]  $S_t$: amount of virgin oil purchased in period $t$.
\end{itemize}
The problem can thus be formulated as the following mixed integer linear program: 
\begingroup
\allowdisplaybreaks
\begin{align}
\rafaelC{[\ IR\ ]} \qquad & z_{IR} = \  \min \ \ \  c\sum_{i \in I}\sum_{\substack{j \in I\\j \neq i}}\sum_{t \in T} d_{ij} X_{ijt} + v\sum_{i\in IC}\sum_{t \in T} X_{0it} + h\sum_{t\in T}I_{0t} + p\sum_{t\in T}S_t  \label{std-obj} & \\
  & \sum_{\substack{j\in I \\ j \neq i }} F_{ijt} -  \sum_{\substack{j\in I \\ j \neq i }} F_{jit} = W_{it}, \qquad  \textrm{for} \ i \in IC, \ t\in T, \label{std-1} \\
&  F_{ijt} \leq (Q-a_{jt})X_{ijt}, \qquad  \textrm{for} \ i \in I, \ j \in I, \ t\in T, \ i\neq j, \label{std-2} \\
&  F_{ijt} \leq Q-W_{jt}, \qquad  \textrm{for} \ i \in I, \ j \in IC, \ t\in T, \ i\neq j, \label{std-3} \\
&  F_{ijt} \geq W_{it} - A_i(1-X_{ijt}), \qquad  \textrm{for} \ i \in IC, \ j \in I, \ t\in T, \ i\neq j, \label{std-4} \\
& \sum_{\substack{j\in I \\ j \neq i }} X_{jit} = Y_{it}, \qquad  \textrm{for} \ i \in IC, \ t\in T, \label{std-5} \\
& \sum_{\substack{j\in I \\ j \neq i }} X_{ijt} = Y_{it}, \qquad  \textrm{for} \ i \in IC, \ t\in T, \label{std-6} \\
& \sum_{i\in IC} X_{i0t} = \sum_{i\in IC} X_{0it}, \qquad  \textrm{for} \ t\in T, \label{std-7} \\
&  W_{it} \leq A_{i}Y_{it}, \qquad  \textrm{for} \ i \in IC, \ t\in T, \label{std-9} \\
&  I_{it} \leq A_{i}(1-Y_{it}), \qquad  \textrm{for} \ i \in IC, \ t\in T, \label{std-10} \\
&  I_{it} = I_{i,t-1} + a_{it}Z_i - W_{it}, \qquad  \textrm{for} \ i \in IC, \ t\in T, \label{std-11} \\
&  I_{i0} = I_{i\tau}, \qquad  \textrm{for} \ i \in I, \label{std-12} \\
&  I_{0t} = I_{0,t-1} + \sum_{i \in IC} W_{it} + S_t - r_t, \qquad  \textrm{for} \ t\in T, \label{std-13} \\
&  Z_{i} \leq \sum_{t\in T} Y_{it}, \qquad  \textrm{for} \ i \in IC, \label{std-14} \\
&  Z_{i} \geq Y_{it}, \qquad  \textrm{for} \ i \in IC, \ t\in T, \label{std-15} \\
&  X_{ijt} \in \{0,1\}, \qquad  \textrm{for} \ i \in I, \ j \in I, \ t\in T, \ i\neq j, \label{std-19} \\
&  Y_{it} \in \{0,1\}, \qquad  \textrm{for} \ i \in IC, \ t\in T, \label{std-20} \\
&  Z_{i} \in \{0,1\}, \qquad  \textrm{for} \ i \in IC, \label{std-21} \\
&  F_{ijt} \geq 0, \qquad  \textrm{for} \ i \in I, \ j \in I, \ t\in T, \ i\neq j, \label{std-22} \\
&  W_{it} \geq 0, \qquad  \textrm{for} \ i \in IC, \ t\in T, \label{std-23} \\
&  I_{it} \geq 0, \qquad  \textrm{for} \ i \in I, \ t\in \{0\} \cup T, \label{std-24} \\
&  S_{t} \geq 0, \qquad  \textrm{for} \ t\in T. \label{std-25}
\end{align}
\endgroup
The objective function \eqref{std-obj} minimizes the total traveling, vehicle usage, inventory and purchase costs. 
Constraints \eqref{std-1} are flow balance constraints for each node $i\in IC$ and guarantee that ${W}_{it}$ liters of oil are collected in the node whenever this amount is nonzero. Constraints \eqref{std-2}-\eqref{std-4} define upper and lower bounds on the flow variables. Constraints \eqref{std-5} and \eqref{std-6} ensure each collection node is visited exactly once whenever it is chosen to be visited. Constraints \eqref{std-7} guarantee that all vehicles leaving the depot in a given period return. 
Constraints \eqref{std-9} imply that oil can only be collected from node $i \in IC$ if it is selected to be visited. Constraints \eqref{std-10} enforce that no inventory remains in collection node $i\in IC$ at the end of a given period if it is visited, i.e., all the oil available in storage must be collected.
Constraints \eqref{std-11} are inventory balance constraints for each collection node $i\in IC$.
Constraints \eqref{std-12} ensure that the initial and ending inventories are the same. Constraints \eqref{std-13} are inventory balance constraints for the depot.
Constraints \eqref{std-14} and \eqref{std-15} set the $Z_i$ variables to zero if node $i\in IC$ is never visited and to one otherwise. Constraints \eqref{std-19}-\eqref{std-25} define integrality and nonnegativity requirements on the variables.

\section{MIP-based heuristic}
\label{sec:mipheuristic}

In this section, we describe a MIP-based heuristic for the selective and periodic inventory routing problem. 
The approach consists of three main steps, which are: (a) solving a relaxation without routing of the problem, (b) constructing routes based on the solution of the relaxation, and (c) improving the routes for each period using a MIP-search. The approach follows a similar structure to the method proposed in \citeA{MonGamGen19} for a related problem, but it differentiates greatly in the way the routes are built and in the improvement of these constructed routes.  
In what follows, Section~\ref{sec:relaxation} describes the relaxation without routing. Section~\ref{sec:routesbuilding} demonstrates how the routes are constructed based on the solution of the relaxation, using heuristics for the traveling salesman problem (TSP)~\cite{AppBixChvCoo06}. Section~\ref{sec:mipsearch} explains the MIP-search for routing improvements, which is based on an exact approach for the capacitated vehicle routing problem (CVRP)~\cite{Lap09}. 

\subsection{Relaxation without routing}\label{sec:relaxation}

The relaxation without routing described in this section is very similar to the one proposed in \citeA{AksKaySibOzg14}, where the authors have shown that such relaxation could provide very strong bounds in short computational times. Basically, the routing related variables and constraints are removed from formulation IR and are replaced by new variables and constraints which guarantee a lower bound on the total routing cost. 
Define the additional decision variables
\begin{itemize}
    \item[--] $V_{t}$: number of vehicles needed to transport the collected oil in period $t\in T$; 
    \item[--] $R_t$: slack capacity to fulfill vehicle with fractional usage in period $t\in T$.
\end{itemize}
The relaxation without routing can thus be formulated as
\begingroup
\allowdisplaybreaks
\begin{align}
\rafaelC{[\ IRR \ ]} \quad & z_{IRR} = \  \min \  c \left[ \sum_{i\in IC} (d_{0i} + d_{i0})\frac{W_{it}}{Q} + \left(\min_{i\in IC}\{d_{0i}\} + \min_{i\in IC}\{d_{i0}\} \right)R_t \right] + v \sum_{t\in T} V_t + h\sum_{t\in T}I_{0t} + p\sum_{t\in T}S_t \label{relax-obj} & \\
&  \rafaelC{\eqref{std-9}-\eqref{std-15}, \eqref{std-20}-\eqref{std-21}, \eqref{std-23}-\eqref{std-25}} \nonumber \\ 
  &  V_{t} = \sum_{i\in IC} \frac{W_{it}}{Q} + R_t, \qquad  \textrm{for} \ t\in T, \label{relax-1} \\
&  W_{it} \leq QY_{it}, \qquad  \textrm{for} \ i\in IC,\ t\in T, \label{relax-8} \\
&  V_{t} \in \rafaelC{\mathbb{Z}_{\geq 0}}, \qquad  \textrm{for} \ t\in T, \label{relax-23a} \\
&  R_{t} \in [0,1], \qquad  \textrm{for} \ t\in T. \label{relax-26}
\end{align}
\endgroup
The objective function \eqref{relax-obj} minimizes a lower bound on the routing cost plus the inventory and purchase costs.
The lower bound can be interpreted as follows. The first term implies that the fraction of the vehicle used to collect waste oil from a given collection node must at least travel from the depot to that node and return back. The second term ensures that the unused capacity of the vehicle must at least travel to and return from the collection nodes which are the closest to and from the depot.
Constraints \eqref{relax-1} set the number of vehicles in each period to the sum of fractional vehicles plus the slack capacity. Constraints \eqref{relax-8} ensure that the total amount collected in a node in a given period does not exceed the capacity of the vehicle. 
\rafaelC{
Constraints~\eqref{relax-23a} define the  integrality and nonnegativity requirements on the $V$ variables while constraints~\eqref{relax-26} ensure the $R$ variables lie in $[0,1]$.
}

Additionally, although not needed in the formulation, we observed that inequalities 
\begin{equation}
\sum_{i\in IC} W_{it} \leq QV_{t}, \qquad  \textrm{for}\ t\in T, \label{relax-extra}
\end{equation}
could help reducing the time to solve several relaxations to optimality using the MIP solver.

\subsection{Routes construction}\label{sec:routesbuilding}

The proposed heuristic to construct routes attempts to build a solution with as few vehicles as possible, starting from a partial solution obtained from solving the relaxation without routing (IRR) described in section \ref{sec:relaxation}.
\rafaelC{Let $\mathcal{S}$ be a feasible/optimal solution to IRR and define $\mathcal{S}.\hat{W}$ and $\mathcal{S}.\hat{Y}$ to represent, respectively, the particular values assigned to the $W$ and $Y$ variables in $\mathcal{S}$. A partial solution defines whether each node $i\in IC$ should be visited in each time period $t\in T$, given by $\mathcal{S}.\hat{Y}_{it}$, as well as the amount of waste oil to be collected, given by $\mathcal{S}.\hat{W}_{it}$.}
The main idea \rafaelC{of the approach} is that the collection nodes to be visited are partitioned into subsets not exceeding the capacity of the vehicles, which are later routed using greedy heuristics for the traveling salesman problem (TSP)~\cite{AppBixChvCoo06}. 
The pseudocode of the approach is described in Algorithm~\ref{alg:buildroutes} and takes as input a partial solution $\mathcal{S}$ obtained from a feasible solution to the relaxation without routing.

Algorithm~\ref{alg:buildroutes} builds routes for each time period independently in the foreach loop of lines~\ref{BR-1}-\ref{BR-1b}. For each time period $t\in T$, the approach firstly initializes the set \rafaelC{$\mathbfcal{B}$} of subsets of collection nodes in line~\ref{BR-2} and the collection nodes to be visited $\mathcal{I}$ in line~\ref{BR-3}.

\rafaelC{In what follows, a heuristic strategy is applied which partitions the collection nodes to be visited seeking for the minimization of the number of subsets in the partition. The main idea of the strategy is to follow a two-stage technique. In the first stage, the algorithm attempts to perform the task heuristically by solving a series of knapsack problems, what can be done very effectively using dynamic programming. In case a lower bound on the number of needed vehicles is achieved, we are done. Otherwise, in the second stage, the task is performed exactly by solving a bin packing problem using an integer programming formulation for which the knapsack based solution obtained in the first stage is given as initial solution.}
The while loop of lines~\ref{BR-4}-\ref{BR-4b} heuristically partitions the collection nodes to be visited in an attempt to minimize the number of subsets in the partition.
While there are collection nodes which are not selected for a subset (line~\ref{BR-4}), a new subset is constructed with those maximizing the usage of the vehicle capacity by solving a knapsack problem. Consider $\alpha_{i}$ to be a binary variable representing whether collection node $i$ is selected. This knapsack problem is formulated as
\begin{align}
\rafaelC{[\ KN (\mathcal{I}, \rafaelC{\mathcal{S}.}\hat{W}) \ ]} \qquad & z_{KN(\mathcal{I}, \rafaelC{\mathcal{S}.}\hat{W})} = \  \max \ \ \  \sum_{i \in \mathcal{I}}\rafaelC{\mathcal{S}.}\hat{W}_{it} \alpha_{i}  \label{knap-obj} & \\
 &  \sum_{i \in \mathcal{I}}\rafaelC{\mathcal{S}.}\hat{W}_{it} \alpha_{i} \leq Q, \label{knap-1} \\
&  \alpha_{i} \in \{0,1\}, \qquad  \textrm{for} \ i \in \mathcal{I}. \label{knap-2} 
\end{align}
The objective function \eqref{knap-obj} maximizes the volume of the selected collection nodes. Constraints \eqref{knap-1} guarantee that the capacity of the vehicle is not exceeded while constraints \eqref{knap-2} ensure the integrality of the variables.
A subset \rafaelC{${B}$} is built by solving this knapsack problem in line \ref{BR-5} using dynamic programming by the auxiliary procedure DP-KNAPSACK ($\mathcal{I},\hat{W}$). Subset $\rafaelC{B}$ is thus inserted into $\mathbfcal{B}$ in line~\ref{BR-6} and the collection nodes in $\rafaelC{B}$ are removed from $\mathcal{I}$ in line~\ref{BR-7}.

Next, in case the partition $\mathbfcal{B}$ contains more subsets (or parts) than a certain lower bound (line~\ref{BR-8}), the approach tries to further reduce the number of subsets by exactly solving a bin packing problem (lines~\ref{BR-9a}-\ref{BR-9}). Considering binary variables $\gamma_{b}$ to define whether bin $b$ is used and $\beta_{ib}$ to indicate whether collection node $i$ is put into bin $b$, this bin packing problem is formulated as
\begin{align}
\rafaelC{[\ BP(\mathcal{I}, \rafaelC{\mathcal{S}.}\hat{W}) \ ]} \qquad & z_{BP(\mathcal{I}, \rafaelC{\mathcal{S}.}\hat{W})} = \  \min \ \ \   \sum_{b \in \rafaelC{\mathbfcal{B}}} \gamma_{b}   \label{binpack-obj} & \\
 &  \sum_{i \in \mathcal{I}}\rafaelC{\mathcal{S}.}\hat{W}_{it} \beta_{ib} \leq Q \gamma_b, \label{binpack-1} \\
& \gamma_{b} \geq \gamma_{b'}, \qquad \textrm{for} \ b<b',  \label{binpack-2} \\
&  \beta_{ib} \in \{0,1\}, \qquad  \textrm{for} \ i \in \mathcal{I}, \ b\in \rafaelC{\mathbfcal{B}}, \label{binpack-3} \\
&  \gamma_{b} \in \{0,1\}, \qquad  \textrm{for} \ b\in \rafaelC{\mathbfcal{B}}. \label{binpack-4}
\end{align}
The objective function \eqref{binpack-obj} minimizes the number of used bins. Constraints \eqref{binpack-1} ensure the capacity of each bin is not exceeded. Constraints \eqref{binpack-2} are symmetry breaking constraints. Finally, constraints \eqref{binpack-3} and \eqref{binpack-4} enforce the integrality requirements on the variables.

In the following, routes are built for each subset individually by heuristically tackling the corresponding traveling salesman problem (lines \ref{BR-10}-\ref{BR-11}). This is performed by using fast heuristics for the TSP, namely, nearest neighbor and farthest insertion. These two greedy heuristics are executed and the best obtained tour amongst the two of them is selected to compose the solution.

Finally, the complete solution is returned in line \ref{BR-12}.

\begin{algorithm}[!ht]
\caption {ROUTES-CONSTRUCTION ($\mathcal{S}$)}
\label{alg:buildroutes}
    
   \ForEach{$t \in T$}{ \label{BR-1}
        $\mathbfcal{B} \leftarrow \varnothing$ \;\label{BR-2}
        $\mathcal{I} \leftarrow \{i \in IC \ | \ \mathcal{S}.\hat{Y}_{it} = 1\} $ \;\label{BR-3}
        \While{$\mathcal{I} \neq \varnothing$}{\label{BR-4}
            Construct subset $\rafaelC{B}$ by solving DP-KNAPSACK ($\mathcal{I}$ , $\mathcal{S}.\hat{W}$)\; \label{BR-5}
            Insert subset $\rafaelC{B}$ into $\mathbfcal{B}$\;\label{BR-6}
            $\mathcal{I} \leftarrow \mathcal{I} \setminus \rafaelC{B}$\;\label{BR-7}
        }\label{BR-4b}
        \If{$|\mathbfcal{B}| > \left\lceil \frac{\sum_{i \in IC} \mathcal{S}.\hat{W}_{it}}{Q}\right\rceil$}{\label{BR-8}
                $\mathcal{I} \leftarrow \{i \in IC \ | \ \mathcal{S}.\hat{Y}_{it} = 1\}$\;\label{BR-9a}
                $\mathbfcal{B} \leftarrow$ \rafaelC{updated} partition constructed from the solution of the bin packing BP ($\mathcal{I}$, $\mathcal{S}.\hat{W}$) solved using a MIP solver \rafaelC{with the original $\mathbfcal{B}$ offered as initial solution}\;\label{BR-9}
           }
       \ForEach{$\rafaelC{B} \in \mathbfcal{B}$}{\label{BR-10}
            Build route for $\mathcal{S}$ in period $t$ visiting collections in $\rafaelC{B}$ using a greedy heuristic for the TSP\;\label{BR-11}
       }\label{BR-10b}
    }\label{BR-1b}
    \Return complete feasible solution $\mathcal{S}$\;\label{BR-12}
\end{algorithm}

\subsection{MIP-search}\label{sec:mipsearch}

The proposed MIP-search is \rafaelC{used as} a postprocessing \rafaelC{local search} procedure which tries to improve \rafaelC{an initial} solution \rafaelC{received as input} by refining the routes traveling to the collection nodes to be visited in a given period.
\rafaelC{The approach receives as input a solution $\mathcal{S}$, the sets of collection nodes to be visited in each period $\mathcal{I}_t$, and the lower bound $lb_t$ on the number of vehicles to be used in each period.
}
This is achieved by solving the associated capacitated vehicle routing problems (CVRPs) for each period \rafaelC{such that the routes in solution $\mathcal{S}$ are offered to the solver as an initial feasible solution (warm start).} In order to achieve this goal, a one-commodity flow mixed integer programming formulation~\cite{GavGra78} is considered for each of these CVRPs. Variables $x_{ij}$ and $f_{ij}$, are the single period counterparts of variables $X_{ijt}$ and $F_{ijt}$. Such CVRP is thus formulated as 
\begin{align}
\rafaelC{[\ CVRP(t,\mathcal{I}_t,\rafaelC{\mathcal{S}.}\hat{W},lb_t) \ ]} \qquad & z_{CVRP(t,\mathcal{I}_t, \rafaelC{\mathcal{S}.}\hat{W}, lb_t)} = \  \min \ \ \  c\sum_{i \in \rafaelC{\mathcal{I}_t \cup \{0\}}}\sum_{\substack{j \in \rafaelC{\mathcal{I}_t \cup \{0\}} \\j \neq i}} d_{ij} x_{ij} + v\sum_{i\in \mathcal{I}_t} x_{0i}  \label{cvrp-obj} & \\
 & \sum_{\substack{j\in \rafaelC{\mathcal{I}_t \cup \{0\}} \\ j \neq i }} f_{ij} -  \sum_{\substack{j\in \rafaelC{\mathcal{I}_t \cup \{0\}} \\ j \neq i }} f_{ji} = \rafaelC{\mathcal{S}.}\hat{W}_{it}, \qquad  \textrm{for} \ i \in \mathcal{I}_t, \label{cvrp-1} \\
&  f_{ij} \leq (Q- \rafaelC{\mathcal{S}.}\hat{W}_{jt})x_{ij}, \qquad  \textrm{for} \ i \in \rafaelC{\mathcal{I}_t \cup \{0\}}, \ j \in \rafaelC{\mathcal{I}_t \cup \{0\}}, \ i\neq j, \label{cvrp-2} \\
&  f_{ij} \geq \rafaelC{\mathcal{S}.}\hat{W}_{it}x_{ij}, \qquad  \textrm{for} \ i \in \mathcal{I}_t, \ j \in \rafaelC{\mathcal{I}_t \cup \{0\}}, \ i\neq j, \label{cvrp-3} \\
& \sum_{\substack{j\in I \\ j \neq i }} x_{ji} = 1, \qquad  \textrm{for} \ i \in \mathcal{I}_t, \label{cvrp-4} \\
& \sum_{\substack{j\in I \\ j \neq i }} x_{ij} = 1, \qquad  \textrm{for} \ i \in \mathcal{I}_t, \label{cvrp-5} \\
& \sum_{i\in \mathcal{I}_t} x_{i0} = \sum_{i\in \mathcal{I}_t} x_{0i}, \label{cvrp-6} \\
& \sum_{i\in \mathcal{I}_t} x_{0i} \geq lb_t, \label{cvrp-lb} \\
&  x_{ij} \in \{0,1\}, \qquad  \textrm{for} \ i \in \mathcal{I}_t\cup \{0\}, \ j \in \mathcal{I}_t\cup \{0\}, \ i\neq j, \label{cvrp-7} \\
&  f_{ij} \geq 0, \qquad  \textrm{for} \ i \in \mathcal{I}_t\cup \{0\}, \ j \in \mathcal{I}_t\cup \{0\}, \ i\neq j. \label{cvrp-8}
\end{align}
The objective function \eqref{cvrp-obj} minimizes the total traveling and vehicle usage costs. Constraints \eqref{cvrp-1} are balance constraints for node $i$ and guarantee $\hat{W}_{it}$ liters of oil are collected in node $i$. Constraints \eqref{cvrp-2} and \eqref{cvrp-3} link the $f$ and $x$ variables, and also define upper and lower bounds on the flow variables. Constraints \eqref{cvrp-4} and \eqref{cvrp-5} ensure each collection point is visited exactly once. Constraints \eqref{cvrp-6} gurantee that all vehicles leaving the depot must return. Constraints~\eqref{cvrp-lb} ensures the number of used vehicles is greater or equal than a lower bound $lb_t$, which is defined by the feasible solution obtained with Algorithm~\ref{alg:buildroutes}. Constraints \eqref{cvrp-7} and \eqref{cvrp-8} impose integrality and nonnegativity requirements on the variables, respectively.

\rafaelC{We remark, though, that the MIP-search procedure described in this section could also be used independently from the initial routes offered by ROUTES-CONSTRUCTION. The drawback of this usage would be the potential increase in running times, as the MIP solver would have to optimize the formulation from scratch without an initial solution. }

\subsection{The complete MIP-based heuristic}

The pseudocode of the MIP-based heuristic is described in Algorithm~\ref{alg:mipheuristic}. Firstly, the relaxation without routing is solved in line~\ref{mipheur-1} using a MIP solver to obtain partial solutions, i.e., without determined routes. In this step, all the solutions encountered by the solver during the search process are stored, and an elite set of solutions containing those which are not more than $\delta \%$ worse than the best one is returned at the end of the execution. The idea behind selecting an elite set of solutions and not simply the best one is to permit variability in the complete constructed solutions. 
After that, routes are constructed for each of these elite partial solutions in lines \ref{mipheur-2}-\ref{mipheur-3}. 
In the following, a set of elite complete solutions is determined with the $k$ best complete solutions (line \ref{mipheur-4}), which are further improved by the MIP-search in lines \ref{mipheur-5}-\ref{mipheur-6}.
The best encountered solution is returned in line~\ref{mipheur-7}.

\begin{algorithm}[H]
\caption {MIP-BASED-HEURISTIC}
\label{alg:mipheuristic}
        ${\rafaelC{\mathbf{\Gamma}}}\leftarrow$ subset of elite partial solutions obtained by solving IRR using a MIP solver\;\label{mipheur-1}
        \ForEach{$ \mathcal{S} \in \rafaelC{\mathbf{\Gamma}}$}{\label{mipheur-2}
           $ \mathcal{S} \leftarrow$ complete solution obtained after executing ROUTES-CONSTRUCTION($\mathcal{S}$)\; \label{mipheur-3}
        }
        $\rafaelC{\mathbf{\Gamma}}' \leftarrow $ subset of elite complete solutions in $\rafaelC{\mathbf{\Gamma}}$\;\label{mipheur-4}
        \ForEach{$\mathcal{S} \in \rafaelC{\mathbf{\Gamma}}' $}{\label{mipheur-5}
            $\mathcal{S} \leftarrow$ improved solution after executing the MIP-search over $\mathcal{S}$ for each period $t\in T$\;\label{mipheur-6}
        }
    \Return $\argmin_{\mathcal{S} \in \rafaelC{\mathbf{\Gamma}}'}\{ z_{IR}(\mathcal{S}) \}$\;\label{mipheur-7}
\end{algorithm}

\rafaelC{
We remark that the proposed MIP-based heuristic provides an \textit{a posteriori} quality guarantee. To see that, define $\underline{z}_{IRR}$ to be the best lower bound at the end of the execution of the MIP solver for solving IRR in line~\ref{mipheur-1}, which is consequently a valid lower bound for IR. Additionally, let $\overline{z}_{IR}$ be the value of the best solution returned in line~\ref{mipheur-7}. Thus, at the end of the execution of Algorithm~\ref{alg:mipheuristic}, it is known that the best encountered solution is within $100 \times \frac{ \overline{z}_{IR} - \underline{z}_{IRR}}{\overline{z}_{IR}}\%$ of optimality.
}

\section{Computational experiments}
\label{sec:experiments}

This section reports the computational experiments conducted to assess the performance of the proposed MIP-based heuristic.
All computational experiments were carried out on a standard Dell Inspiron 15 laptop, running under Ubuntu GNU/Linux, with an Intel(R) Core(TM) i7-7500U CPU @ 2.70GHz processor and 8GB of RAM. The algorithms were coded in Julia v1.4.0, using JuMP v0.18.6. The formulations were solved using Gurobi 8.0.1 with the standard configurations, except the relative optimality tolerance gap which was set to $10^{-6}$. A time limit of 60 seconds was imposed for every execution of the MIP solver.

\subsection{Benchmark instances}\label{sec:instancias}

The considered benchmark instance set is composed of three benchmark sets used in the literature. In all the considered instances, the cyclic planning horizon consists of seven days.

The first benchmark set was introduced in \citeA{AksKaySalAkc12}, containing 36 small instances with 25 collection nodes. There are two types of vehicles, namely: Fiat Doblo Cargo Maxi (Dob) and Fiat Fiorino Cargo (Fio).
Dob has a capacity ($Q$) of 920 liters, operating cost ($v$) of 110 monetary units per day and per km traveling cost ($c$) of 0.22 monetary units. Fio has a capacity ($Q$) of 550 liters, operating cost ($v$) of 90 monetary units per day and per km traveling cost ($c$) of 0.19 monetary units. 
The waste vegetable oil inventory cost ($h$) is 0.02 monetary units per liter per day. Virgin vegetable oil prices ($p$) lie in $\{0.25, 0.50, 1.25\}$ monetary units per liter. Daily accumulation rates are either low (30 liters per day) or high (60 liters per day). Three levels of oil requirements are established for each value of accumulation rate. For accumulation rate 30, waste requirements can be low (600 liters), medium , or high (900 liters). For accumulation rate 60, waste requirements can be low (1200 liters), medium (1500 liters), or high (1800 liters).

The second benchmark set was proposed in \citeA{AksKaySibOzg14}, containing 54 small and medium instances based on a real-world case study with 20, 25, 30, 35, 40, 50, 60, 80, and 100 collection nodes. There is a single type of vehicle with a capacity ($Q$) of 550 liters. The operating cost of a vehicle ($v$) is 90 monetary units per day while the per km traveling cost ($c$) corresponds to 0.24 monetary units. The waste vegetable oil inventory cost ($h$) is 0.02 monetary units per liter per day, while the virgin oil purchasing price ($p$) lies in $\{2.5,3.5\}$ monetary units per liter. Three levels of oil requirements are established as low, medium, and high.
The instances are identified as ($P_1$n-$P_2$r-$P_3$p), where $P_1$ denotes the number of collection nodes, $P_2$ gives the per day oil requirement and $P_3$ represents the virgin oil purchasing price. 

The third benchmark set was introduced in \citeA{CarGonTreGar19}, containing 24 large instances with 120, 160, 200, and 300 collection nodes. These instances were generated based on those proposed in \cite{AksKaySibOzg14}.

\subsection{Considered approaches and settings}
The following approaches are considered in the comparisons:
\begin{itemize}
    \item[--] MIP-based heuristic in which the MIP-search is not executed (MH)\rafaelC{, i.e., Algorithm~\ref{alg:mipheuristic} (lines \ref{mipheur-1}-\ref{mipheur-4} and \ref{mipheur-7})};
    \item[--] MIP-based heuristic with MIP-search (MH$^{+}$)\rafaelC{, i.e., Algorithm~\ref{alg:mipheuristic} (lines \ref{mipheur-1}-\ref{mipheur-7})};
    \item[--] Adaptive large neighborhood search (ALNS)\cite{AksKaySibOzg14};
    \item[--] Reduce and optimize approach (ROA)\cite{CarGonTreGar19}.
\end{itemize}

For MH and MH$^{+}$, the subset of elite partial solutions obtained by solving IRR (Algorithm~\ref{alg:mipheuristic}, line~\ref{mipheur-1}) is composed of those whose values lie within $\delta = 5\%$ of the best encountered one. For MH$^{+}$, the subset of elite solutions in $\mathcal{S}$ (Algorithm~\ref{alg:mipheuristic}, line~\ref{mipheur-4}) is composed of \rafaelC{simply the best solution as the goal} was to have an approach which could deliver good quality results within low computational times.

\rafaelC{
The values presented for ALNS and ROA correspond to those reported in their corresponding works. 
As the running time of each execution of ROA was fixed in 36 minutes (2160 seconds), we do not report this information in our tables.
}

\subsection{Results}

Table \ref{tab:tab1} presents the results obtained with MH, MH$^{+}$, and ROA for the instances from the first benchmark set. The first column identifies the instance. The next six columns report, for MH and MH$^{+}$, the best solution value, the running time, and the optimality gap obtained with each of the approaches, calculated as $100 \times \rafaelC{(upperbound - lowerbound)/upperbound}$, where \rafaelC{$upperbound$} is the encountered heuristic value and \rafaelC{$lowerbound$} is the bound obtained by the relaxation without routing. The last two columns give the best solution value and the optimality gap obtained with ROA. We remark that \citeA{AksKaySibOzg14} did not report results for ALNS on this instance set. The results show that MH could obtain a low average gap (3.3\%) within a very short average computational time (19.7 seconds). MH$^{+}$ could further improve the results obtained by MH, while keeping the average running time below 40.0 seconds. Note that although ROA obtains a lower average gap (1.6\%), these results were obtained using a much larger computational time of 2160.0 seconds on average, which represents nearly 110$\times$ the average running time of MH. Furthermore, new best results were obtained for 16.7\% of the instances.

Tables \ref{tab:tab2} and \ref{tab:tab3} show the results obtained with MH, MH$^{+}$, ALNS and ROA for the second benchmark set. 
The first column identifies the instance. The next nine columns report, for MH, MH$^{+}$ and ALNS, the best solution value, the running time and the optimality gap obtained with each of the approaches. The last two columns give the best solution value and the optimality gap obtained with ROA. We remark that there are two lines with a problem in the values reported in \cite{AksKaySibOzg14} for ALNS, which are marked with an '$^*$', as their reported values are below the bound provided by the relaxation without routing. 
These tables show the effectiveness of our approach. 
Table \ref{tab:tab2} shows that, for the small instances, MH found an average gap of 3.0\% within an average time of 16.9 seconds, and MH$^{+}$ could further reduce this average gap to 2.7\% within an average time of 25.9 seconds. Both results improve over that obtained by ALNS, and are very close to the one obtained by ROA which uses a much larger computational time of 2160 seconds. Besides, MH and MH$^{+}$ improved the result of ALNS for, respectively, 53.3\% and 66.7\% of the instances. These improvements over the results of ROA were, respectively, 23.3\% and 26.7\%.
Table \ref{tab:tab3} shows that the improvements achieved by MH and MH$^{+}$ over ALNS are even more evident, as better solutions were obtained for, respectively, 75.0\% and 91.7\% of the instances. It is remarkable that the average gap obtained by MH$^{+}$ is much lower than that of ALNS and is very close to that of ROA, but spending just a fraction of the time on average.
Furthermore, new best solutions were encountered for 25.9\% of these instances.

Table \ref{tab:tab4} presents the results obtained with MH, MH$^{+}$ and ROA for the instances from the third benchmark set. The columns are the same as in Table~\ref{tab:tab1}. The results show the impressive results of our heuristics, as they could improve the best known result for all the instances. MH was able to obtain an average gap of 1.6\% within 90.9 seconds on average, while MH$^{+}$ could reduce this gap to 1.4\% within 469.6 seconds on average. These average gaps are much lower than the 5.9\% average gap obtained by ROA using an average time of 2160  seconds.

\begin{table}[H]
\small
\centering
\caption{Results for MH, MH$^{+}$ and ROA on instances from the first benchmark set.}
\begin{tabular}{ l | c c c | c c c | c c }
\hline
\multicolumn{1}{l|}{Instances} &  \multicolumn{3}{c|}{MH} &  \multicolumn{3}{c|}{MH$^{+}$} &  \multicolumn{2}{c}{ROA}\\
 & $z_{MH}$& time$_{MH}$ & gap$_{MH}$ & $z_{MH^{+}}$& time$_{MH^{+}}$ & gap$_{MH^{+}}$ & $z_{ROA}$ & gap$_{ROA}$\\ \hline
Fio-30acc-LOW-025  & 867.4 & 14.6 & 1.6 & 866.3 & 15.3 & 1.5 & \textbf{865.7} & 1.5\\
Fio-30acc-MED-025  & 1099.6 & 17.2 & 1.5 & 1098.6 & 18.0 & 1.4 & \textbf{1094.6} & 1.1\\
Fio-30acc-HIGH-025  & 1348.3 & 15.3 & 2.6 & 1345.0 & 22.3 & 2.4 & \textbf{1331.5} & 1.4\\
Fio-60acc-LOW-025  & 1851.7 & 21.8 & 7.8 & 1847.6 & 33.8 & 7.6 & \textbf{1748.9} & 2.4\\
Fio-60acc-MED-025  & 2273.8 & 16.1 & 6.1 & 2258.8 & 31.0 & 5.4 & \textbf{2233.8} & 4.4\\
Fio-60acc-HIGH-025  & 2874.7 & 13.7 & 8.5 & 2869.2 & 17.2 & 8.3 & \textbf{2713.3} & 3.0\\
Dob-30acc-LOW-025  & 701.1 & 15.3 & 1.1 & \textbf{699.2} & 16.1 & 0.8 & 701.5 & 1.2\\
Dob-30acc-MED-025  & 856.7 & 13.9 & 2.1 & 855.1 & 14.8 & 1.9 & \textbf{847.4} & 1.1\\
Dob-30acc-HIGH-025  & 1024.0 & 14.5 & 2.1 & 1022.5 & 15.7 & 1.9 & \textbf{1019.0} & 1.6\\
Dob-60acc-LOW-025  & 1299.7 & 52.8 & 2.1 & \textbf{1297.9} & 53.1 & 1.9 & 1300.2 & 2.1\\
Dob-60acc-MED-025  & 1636.8 & 15.5 & 2.3 & \textbf{1632.5} & 19.7 & 2.0 & 1647.1 & 2.9\\
Dob-60acc-HIGH-025  & 2144.5 & 17.6 & 2.5 & 2140.0 & 21.5 & 2.3 & \textbf{2133.2} & 2.0\\
Fio-30acc-LOW-050  & 907.4 & 31.9 & 1.2 & 907.0 & 32.2 & 1.1 & \textbf{902.4} & 0.6\\
Fio-30acc-MED-050  & 1140.0 & 15.7 & 1.8 & 1135.1 & 16.8 & 1.3 & \textbf{1127.9} & 0.7\\
Fio-30acc-HIGH-050  & 1539.3 & 14.8 & 1.7 & \textbf{1536.6} & 16.3 & 1.6 & 1549.2 & 2.4\\
Fio-60acc-LOW-050  & 1996.3 & 68.5 & 12.6 & 1985.7 & 134.2 & 12.1 & \textbf{1782.6} & 2.1\\
Fio-60acc-MED-050  & 2336.1 & 15.1 & 8.0 & 2328.1 & 146.4 & 7.7 & \textbf{2236.2} & 3.9\\
Fio-60acc-HIGH-050  & 3253.9 & 13.4 & 2.6 & \textbf{3247.2} & 239.8 & 2.4 & 3261.8 & 2.8\\
Dob-30acc-LOW-050  & 705.9 & 14.8 & 1.2 & 704.8 & 15.6 & 1.1 & \textbf{701.2} & 0.5\\
Dob-30acc-MED-050  & 853.9 & 15.4 & 1.8 & 853.1 & 16.2 & 1.7 & \textbf{845.5} & 0.8\\
Dob-30acc-HIGH-050  & 1220.3 & 14.0 & 1.8 & 1219.3 & 14.9 & 1.8 & \textbf{1218.0} & 1.6\\
Dob-60acc-LOW-050  & 1320.3 & 37.8 & 1.3 & \textbf{1319.2} & 38.4 & 1.3 & 1369.5 & 4.9\\
Dob-60acc-MED-050  & 1692.7 & 15.9 & 3.5 & 1685.7 & 31.6 & 3.1 & \textbf{1647.7} & 0.9\\
Dob-60acc-HIGH-050  & 2817.7 & 12.8 & 5.9 & 2814.2 & 14.0 & 5.8 & \textbf{2657.0} & 0.2\\
Fio-30acc-LOW-125  & 911.3 & 33.8 & 1.6 & 907.6 & 34.6 & 1.2 & \textbf{902.5} & 0.7\\
Fio-30acc-MED-125  & 1136.6 & 14.5 & 1.5 & 1132.8 & 16.0 & 1.1 & \textbf{1125.8}& 0.5\\
Fio-30acc-HIGH-125  & 2064.7 & 13.0 & 2.3 & 2054.9 & 15.2 & 1.8 & \textbf{2025.4} & 0.4\\
Fio-60acc-LOW-125  & 1897.7 & 19.1 & 6.6 & 1892.7 & 24.8 & 6.4 & \textbf{1781.0} & 0.5\\
Fio-60acc-MED-125  & 2390.5 & 16.1 & 8.4 & 2382.6 & 165.7 & 8.1 & \textbf{2242.5} & 2.4\\
Fio-60acc-HIGH-125  & 4972.4 & 15.5 & 3.8 & 4967.6 & 16.8 & 3.7 & \textbf{4841.2} & 1.2\\
Dob-30acc-LOW-125  & 707.0 & 15.9 & 1.4 & 706.3 & 16.8 & 1.3 & \textbf{703.4} & 0.9\\
Dob-30acc-MED-125  & 853.7 & 15.3 & 1.8 & 852.8 & 16.0 & 1.7 & \textbf{845.6} & 0.9\\
Dob-30acc-HIGH-125  & 1765.1 & 13.1 & 1.9 & 1764.4 & 13.9 & 1.9 & \textbf{1740.3} & 0.6\\
Dob-60acc-LOW-125  & 1367.6 & 27.8 & 1.0 & 1364.9 & 28.7 & 0.8 & \textbf{1361.5} & 0.5\\
Dob-60acc-MED-125  & 1680.8 & 14.3 & 2.7 & 1675.2 & 15.0 & 2.3 & \textbf{1654.4} & 1.1\\
Dob-60acc-HIGH-125  & 4390.6 & 12.8 & 3.7 & 4383.5 & 14.0 & 3.5 & \textbf{4233.2} & 0.1 \\ \hline
Averages  &  & 19.7 & 3.3 &  & 38.1 & 3.1 &  & 1.6 \\ 
$<z_{ROA}$(\%)  & 16.7 &  &  & 16.7 &  &  &  &  \\ \hline
\end{tabular}
\label{tab:tab1}
\end{table}

\begin{landscape}

\begin{table}[H]
\small
\centering
\caption{Results for MH, MH$^{+}$, ALNS and ROA on the small instances from the second benchmark set.}
\begin{tabular}{ l | c c c | c c c | c c c | c c }
\hline
\multicolumn{1}{l|}{Instances} &  \multicolumn{3}{c|}{MH} &  \multicolumn{3}{c|}{MH$^{+}$} & \multicolumn{3}{c|}{ALNS} &  \multicolumn{2}{c}{ROA}\\
  & $z_{MH}$& time$_{MH}$ & gap$_{MH}$ & $z_{MH^{+}}$& time$_{MH^{+}}$ & gap$_{MH^{+}}$ & $z_{ALNS}$  & time$_{ALNS}$  &    gap$_{ALNS}$  & $z_{ROA}$ & gap$_{ROA}$\\ \hline
20n-270r-2.5p  & 474.1 & 13.9 & 0.8 & \textbf{474.0} & 14.6 & 0.7 & 480.1 & 3.2 & 2.0 &  475.1 & 1.0\\
20n-410r-2.5p  & 714.2 & 15.1 & 1.8 & 713.7 & 15.8 & 1.7 & 711.8 & 2.8 & 1.5  & \textbf{711.7} & 1.4\\
20n-540r-2.5p  & 856.7 & 12.7 & 2.7 & 856.7 & 13.5 & 2.7 & \textbf{830.6} & 8.6 & -0.3$^*$ &  841.2 & 0.9\\
20n-270r-3.5p & 474.9 & 13.5 & 0.9 & \textbf{474.6} & 14.3 & 0.9 & 479.9 & 3.9 & 2.0 &  475.1 & 1.0\\
20n-410r-3.5p  & 713.5 & 14.4 & 1.7 & 713.3 & 15.1 & 1.7 & 712.0 & 3.7 & 1.5 &  \textbf{708.8} & 1.0\\
20n-540r-3.5p  & 868.0 & 12.8 & 2.7 & 868.0 & 13.5 & 2.7 & \textbf{830.6} & 9.8 & -1.7$^*$ &  852.3 & 0.9\\
25n-320r-2.5p  & 593.1 & 14.2 & 1.4 & 593.0 & 15.0 & 1.4 & 593.2 & 4.7 & 1.4 &  \textbf{592.2} & 1.2\\
25n-480r-2.5p  & 820.9 & 15.3 & 2.5 & 820.0 & 16.2 & 2.4 & 813.0 & 7.1 & 1.5 &  \textbf{808.8} & 1.0\\
25n-640r-2.5p  & 1103.6 & 15.5 & 4.2 & 1102.0 & 16.4 & 4.0 & 1074.1 & 25.4 & 1.5 &  \textbf{1071.3} & 1.3\\
25n-320r-3.5p  & 595.3 & 14.9 & 1.7 & 594.4 & 15.8 & 1.6 & 594.6 & 5.1 & 1.6 &  \textbf{593.8} & 1.5\\
25n-480r-3.5p  & 817.8 & 14.6 & 2.1 & 815.3 & 15.4 & 1.8 & 820.2 & 6.7 & 2.4 &  \textbf{809.5} & 1.1\\
25n-640r-3.5p  & 1093.0 & 15.2 & 3.0 & 1089.3 & 16.0 & 2.7 & 1083.4 & 35.3 & 2.2 &  \textbf{1073.5} & 1.3\\
30n-420r-2.5p  & 709.6 & 21.4 & 2.1 & \textbf{708.8} & 22.1 & 2.0 & 709.4 & 12.1 & 2.1 &  709.0 & 2.1\\
30n-630r-2.5p  & 998.6 & 16.8 & 3.0 & \textbf{996.1} & 19.0 & 2.8 & 1008.4 & 15.2 & 4.0 &  1050.6 & 7.8\\
30n-840r-2.5p  & 1344.3 & 14.5 & 3.9 & 1341.1 & 15.5 & 3.7 & 1420.4 & 224.4 & 9.0 &  \textbf{1323.4} & 2.4\\
30n-420r-3.5p  & 713.5 & 25.3 & 2.7 & 712.7 & 26.1 & 2.6 & 713.5 & 14.0 & 2.7 &  \textbf{710.3} & 2.2\\
30n-630r-3.5p  & 1009.1 & 18.4 & 3.0 & \textbf{1007.7} & 79.3 & 2.9 & 1066.3 & 14.1 & 8.2 &  1053.5 & 7.1\\
30n-840r-3.5p  & 1353.7 & 17.3 & 4.3 & 1345.8 & 27.8 & 3.8 & 1421.5 & 170.9 & 8.9 &  \textbf{1326.0} & 2.3\\
35n-480r-2.5p  & 812.0 & 15.9 & 3.4 & 810.6 & 16.7 & 3.2 & 809.0 & 6.3 & 3.0 &  \textbf{799.5} & 1.8\\
35n-710r-2.5p  & 1143.3 & 24.5 & 3.6 & \textbf{1139.5} & 65.4 & 3.3 & 1178.3 & 13.7 & 6.4 &  1166.5 & 5.5\\
35n-950r-2.5p  & 1557.8 & 15.3 & 3.4 & 1551.1 & 16.2 & 3.0 & 1554.6 & 96.1 & 3.2 &  \textbf{1520.7} & 1.0\\
35n-480r-3.5p  & 812.6 & 16.9 & 3.4 & 811.9 & 17.8 & 3.3 & 810.8 & 8.7 & 3.2 &  \textbf{796.4} & 1.5\\
35n-710r-3.5p  & 1175.5 & 16.1 & 4.5 & 1164.5 & 90.0 & 3.6 & 1180.5 & 10.9 & 4.9 &  \textbf{1158.7} & 3.1\\
35n-950r-3.5p  & 1565.4 & 18.0 & 3.7 & 1542.6 & 20.2 & 2.3 & 1581.2 & 79.6 & 4.7 &  \textbf{1525.5} & 1.2\\
40n-550r-2.5p  & 819.4 & 16.0 & 3.7 & \textbf{817.4} & 16.9 & 3.5 & 832.4 & 7.5 & 5.2 &  831.8 & 5.1\\
40n-820r-2.5p  & 1310.9 & 26.5 & 3.4 & 1300.7 & 30.8 & 2.6 & 1334.5 & 12.3 & 5.1 &  \textbf{1280.6} & 1.1\\
40n-1090r-2.5p  & 1692.9 & 13.4 & 4.7 & 1672.2 & 63.9 & 3.5 & 1673.4 & 141.1 & 3.5 &  \textbf{1641.6} & 1.7\\
40n-550r-3.5p  & 823.5 & 17.8 & 4.2 & \textbf{820.7} & 18.8 & 3.9 & 838.5 & 9.5 & 5.9 &  844.3 & 6.6\\
40n-820r-3.5p  & 1314.7 & 26.4 & 3.6 & 1305.7 & 30.1 & 3.0 & 1373.5 & 11.8 & 7.8 &  \textbf{1287.7} & 1.6\\
40n-1090r-3.5p  & 1697.5 & 13.7 & 4.8 & 1674.9 & 18.2 & 3.5 & 1671.5 & 239.8 & 3.3 &  \textbf{1650.7} & 2.1\\ \hline
Averages  &  & 16.9 & 3.0 &  & 25.9 & 2.7 &  & 40.1 & 3.6 &   & 2.3\\
$<z_{ALNS}$ (\%)  & 53.3 &  &  & 66.7 &  &  &  &  &  &   & \\
$<z_{ROA}$ (\%)  & 23.3 &  &  & 26.7 &  &  &  &  &  &   & \\ \hline
\end{tabular}
\label{tab:tab2}
\end{table}

\begin{table}[H]
\small
\centering
\caption{Results for MH, MH$^{+}$, ALNS and ROA on the medium instances from the second benchmark set.}
\begin{tabular}{ l | c c c | c c c | c c c | c c }
\hline
\multicolumn{1}{l|}{Instances} &  \multicolumn{3}{c|}{MH} &  \multicolumn{3}{c|}{MH$^{+}$} & \multicolumn{3}{c|}{ALNS} &  \multicolumn{2}{c}{ROA}\\
  & $z_{MH}$& time$_{MH}$ & gap$_{MH}$ & $z_{MH^{+}}$& time$_{MH^{+}}$ & gap$_{MH^{+}}$ & $z_{ALNS}$& time$_{ALNS}$ & gap$_{ALNS}$ & $z_{ROA}$ & gap$_{ROA}$\\ \hline
50n-650r-2.5p  & 1076.9 & 20.2 & 4.1 & 1073.8 & 21.0 & 3.8 & 1069.3 & 10.2 & 3.4 & \textbf{1055.5} & 2.1\\
50n-910r-2.5p  & 1453.7 & 17.1 & 4.6 & 1436.2 & 89.2 & 3.5 & 1447.1 & 19.2 & 4.2 & \textbf{1396.0} & 0.7\\
50n-1300r-2.5p  & 2057.8 & 21.5 & 3.9 & 2032.3 & 72.7 & 2.7 & 2129.9 & 286.4 & 7.1 & \textbf{2011.1} & 1.7\\
50n-650r-3.5p  & 1066.4 & 26.4 & 3.1 & 1062.5 & 26.7 & 2.8 & 1070.2 & 12.4 & 3.5 & \textbf{1054.6} & 2.0\\
50n-910r-3.5p  & 1443.1 & 18.8 & 3.9 & 1431.1 & 29.6 & 3.1 & 1489.6 & 29.7 & 6.9 & \textbf{1405.2} & 1.3\\
50n-1300r-3.5p  & 2059.2 & 21.3 & 4.0 & 2033.1 & 29.1 & 2.7 & 2134.2 & 291.4 & 7.3 & \textbf{2025.9} & 2.4\\
60n-800r-2.5p  & 1290.1 & 29.0 & 3.6 & 1283.1 & 30.1 & 3.1 & 1320.8 & 17.4 & 5.9 & \textbf{1265.0} & 1.7\\
60n-1190r-2.5p  & 1885.6 & 18.0 & 3.3 & 1872.0 & 21.6 & 2.6 & 1936.5 & 48.6 & 5.9 & \textbf{1851.8} & 1.6\\
60n-1590r-2.5p  & 2513.9 & 20.0 & 3.8 & 2490.7 & 27.0 & 2.9 & 2544.9 & 432.4 & 5.0 & \textbf{2460.0} & 1.7\\
60n-800r-3.5p  & 1286.4 & 28.6 & 3.3 & 1281.3 & 30.2 & 3.0 & 1326.5 & 19.3 & 6.3 & \textbf{1264.1} & 1.6\\
60n-1190r-3.5p  & 1910.3 & 16.3 & 4.6 & 1878.6 & 25.3 & 3.0 & 1954.1 & 55.1 & 6.7 & \textbf{1857.3} & 1.8\\
60n-1590r-3.5p  & 2523.9 & 19.7 & 3.9 & 2488.9 & 31.8 & 2.5 & 2562.6 & 395.0 & 5.3 & \textbf{2471.9} & 1.8\\
80n-1070r-2.5p  & 1634.0 & 36.0 & 4.2 & 1622.4 & 38.4 & 3.6 & 1722.8 & 59.5 & 9.2 & \textbf{1599.8} & 2.2\\
80n-1610r-2.5p  & 2461.4 & 28.5 & 3.6 & \textbf{2435.9} & 37.8 & 2.6 & 2535.4 & 173.2 & 6.4 & 2462.7 & 3.7\\
80n-2150r-2.5p  & 3491.9 & 13.4 & 5.5 & 3421.5 & 225.3 & 3.5 & 3473.1 & 979.5 & 5.0 & \textbf{3361.1} & 1.8\\
80n-1070r-3.5p  & 1632.9 & 42.7 & 4.2 & 1617.8 & 45.5 & 3.3 & 1732.8 & 53.0 & 9.7 & \textbf{1598.5} & 2.1\\
80n-1610r-3.5p  & 2457.6 & 24.2 & 3.5 & \textbf{2429.5} & 34.3 & 2.4 & 2581.9 & 198.7 & 8.1 & 2434.5 & 2.6\\
80n-2150r-3.5p  & 3532.9 & 17.2 & 5.5 & 3458.3 & 211.5 & 3.4 & 3512.1 & 1038.3 & 4.9 & \textbf{3419.2} & 2.3\\
100n-1330r-2.5p  & 2024.0 & 44.5 & 5.1 & \textbf{2005.1} & 161.7 & 4.2 & 2083.5 & 89.1 & 7.8 & 2063.3 & 6.9\\
100n-2000r-2.5p  & 3082.0 & 43.8 & 4.4 & \textbf{3048.2} & 358.8 & 3.3 & 3149.7 & 256.3 & 6.4 & 3149.5 & 6.4\\
100n-2670r-2.5p  & 4313.7 & 19.7 & 7.0 & 4297.3 & 440.7 & 6.7 & 4293.2 & 3492.3 & 6.6 & \textbf{4234.8} & 5.3\\
100n-1330r-3.5p  & 2021.0 & 50.5 & 5.0 & \textbf{1994.2} & 190.7 & 3.7 & 2149.1 & 77.4 & 10.7 & 2060.2 & 6.8\\
100n-2000r-3.5p  & 3094.8 & 47.0 & 4.8 & \textbf{3060.7} & 327.5 & 3.7 & 3233.7 & 295.1 & 8.9 & 3164.8 & 6.9\\
100n-2670r-3.5p  & 4335.5 & 20.5 & 6.7 & 4327.4 & 441.6 & 6.5 & 4333.3 & 3294.8 & 6.6 & \textbf{4253.6} & 4.9\\ \hline
Averages  &  & 26.9 & 4.4 &  & 122.8 & 3.4 &  & 484.3 & 6.6 &   & 3.0\\
$<z_{ALNS}$ (\%)  & 75.0 &  &  & 91.7 &  &  &  &  &  &   & \\
$<z_{ROA}$ (\%)  & 20.0 &  &  & 24.0 &  &  &  &  &  &   & \\ \hline
\end{tabular}
\label{tab:tab3}
\end{table}

\end{landscape}

\begin{table}[H]
\footnotesize
\centering
\caption{Results for MH, MH$^{+}$ and ROA on the very large instances from the third benchmark set.}
\begin{tabular}{ l | c c c | c c c | c c }
\hline
\multicolumn{1}{l|}{Instances} &  \multicolumn{3}{c|}{MH} &  \multicolumn{3}{c|}{MH$^{+}$} &  \multicolumn{2}{c}{ROA}\\
  & $z_{MH}$& time$_{MH}$ & gap$_{MH}$ & $z_{MH^{+}}$& time$_{MH^{+}}$ & gap$_{MH^{+}}$ & $z_{ROA}$ & gap$_{ROA}$\\ \hline
120n-1410r-2.5p  & 1686.2 & 73.4 & 1.4 & \textbf{1681.6} & 343.9 & 1.1 & 1766.7 & 5.9\\
120n-2030r-2.5p  & 2424.2 & 48.3 & 1.4 & \textbf{2416.6} & 325.9 & 1.1 & 2565.6 & 6.9\\
120n-2820r-2.5p  & 3404.0 & 82.9 & 2.8 & \textbf{3393.2} & 504.7 & 2.5 & 3530.4 & 6.3\\
120n-1410r-3.5p &  1687.3 & 63.5 & 1.4 & \textbf{1683.5} & 294.7 & 1.2 & 1763.9 & 5.7\\
120n-2030r-3.5p &  2426.0 & 59.9 & 1.5 & \textbf{2420.0} & 434.7 & 1.3 & 2569.1 & 7.0\\
120n-2820r-3.5p  & 3376.2 & 74.3 & 2.0 & \textbf{3373.9} & 495.8 & 1.9 & 3514.7 & 5.9\\
160n-1810r-2.5p  & 2208.2 & 73.2 & 1.9 & \textbf{2205.6} & 493.8 & 1.8 & 2229.3 & 2.8\\
160n-2930r-2.5p  & 3525.7 & 78.8 & 1.5 & \textbf{3514.3} & 401.7 & 1.2 & 3702.8 & 6.3\\
160n-3320r-2.5p  & 4039.6 & 91.9 & 3.0 & \textbf{4028.5} & 513.1 & 2.7 & 4157.9 & 5.7\\
160n-1810r-3.5p  & 2216.0 & 73.2 & 1.3 & \textbf{2209.2} & 277.5 & 1.0 & 2235.8 & 2.2\\
160n-2930r-3.5p  & 3518.1 & 149.1 & 1.3 & \textbf{3509.7} & 495.1 & 1.1 & 3624.4 & 4.2\\
160n-3320r-3.5p  & 3983.8 & 93.2 & 1.6 & \textbf{3977.8} & 460.0 & 1.5 & 4066.1 & 3.6\\
200n-2110r-2.5p  & 2505.7 & 74.5 & 1.3 & \textbf{2498.9} & 451.9 & 1.0 & 2575.1 & 3.9\\
200n-3430r-2.5p  & 4068.5 & 60.7 & 1.2 & \textbf{4064.4} & 481.6 & 1.1 & 4257.9 & 5.6\\
200n-4120r-2.5p  & 4899.6 & 154.3 & 1.3 & \textbf{4891.9} & 575.2 & 1.2 & 5168.3 & 6.5\\
200n-2110r-3.5p  & 2507.8 & 75.1 & 1.3 & \textbf{2505.0} & 496.0 & 1.2 & 2575.3 & 3.9\\
200n-3430r-3.5p  & 4069.5 & 68.3 & 1.2 & \textbf{4061.4} & 489.9 & 1.0 & 4250.4 & 5.4\\
200n-4120r-3.5p 8 & 4907.2 & 93.5 & 1.5 & \textbf{4901.4} & 514.4 & 1.4 & 5087.7 & 5.0\\
300n-2510r-2.5p 9 & 3018.1 & 76.1 & 1.3 & \textbf{3013.5} & 497.1 & 1.1 & 3207.5 & 7.1\\
300n-4030r-2.5p  & 5014.6 & 152.2 & 4.0 & \textbf{5004.3} & 574.0 & 3.8 & 5262.3 & 8.5\\
300n-4520r-2.5p  & 5433.6 & 96.2 & 1.0 & \textbf{5430.8} & 517.5 & 1.0 & 5865.7 & 8.3\\
300n-2510r-3.5p  & 3025.2 & 74.9 & 1.5 & \textbf{3025.0} & 495.8 & 1.5 & 3197.9 & 6.8\\
300n-4030r-3.5p  & 4861.7 & 134.8 & 0.9 & \textbf{4856.5} & 555.7 & 0.8 & 5258.7 & 8.4\\
300n-4520r-3.5p  & 5443.0 & 159.1 & 1.2 & \textbf{5437.6} & 580.0 & 1.1 & 5973.3 & 10.0\\ \hline
Averages  &  & 90.9 & 1.6 &  & 469.6 & 1.4 &  & 5.9 \\
$<z_{ROA}$ (\%)  & 100.0 &  &  & 100.0 &  &  &  &  \\ \hline
\end{tabular}
\label{tab:tab4}
\end{table}


\rafaelC{
Tables~\ref{tab:statsfio}-\ref{tab:statslarge} summarize the executions of ROUTES-CONSTRUCTION (Algorithm~\ref{alg:buildroutes}) and the characteristics of the best obtained solutions using MH$^+$ for each instance.
We observed that the values which will be presented for the best solutions matched for MH and MH$^+$ in nearly all instances and, thus, the improvements are more likely to be related to reorganizations of the routes. Therefore, only the solutions obtained by MH$^+$ are summarized.
In each of these tables, the first column gives the instance. Columns 2-4 present information about the execution of procedure ROUTES-CONSTRUCTION, which are the number of times the heuristic knapsack based partitioning was performed (\#knpart), the number of times the bin packing based partitioning was performed (\#bppart), and the number of times in which the latter improved the former (\#bpimpr).
Next, columns 5-10 outline the best encountered solutions, with columns 5-7 presenting information about the number of vehicles per period (minimum, average, and maximum) and columns 8-10 displaying information about the number of collections per vehicle (minimum, average, and maximum).
The last line of each table presents the averages over the columns.
We remark that whenever all collections to be performed based on the partial solution $\mathcal{S}$ could fit in a single vehicle, \#knpart is given as 0.
}

\rafaelC{
It can be observed from the tables that for several of the small instances, the approach used at most one vehicle per period. Besides, for some of them, collection did not happen in every period, what can be seen by the value 0 in column min(veh).
}
\rafaelC{A noteworthy observation is the fact that the average values in the columns related to ROUTES-CONSTRUCTION executions in the tables consistently increased as the sizes of the instances grew larger.}
\rafaelC{It is also interesting to observe that, for small and medium instances, the heuristic dynamic programming based partitioning was very effective and could rarely be improved by the bin packing based partitioning. On the other hand, the bin packing based partitioning allowed improvements for several executions of ROUTES-CONSTRUCTION for the large instances, showing itself very valuable.
}
\rafaelC{
Furthermore, we can see that for the large instances the numbers of used vehicles per period are considerably larger than those for the small and medium ones. This is a possible reason why the gaps obtained for the large instances are very low, as the costs of the solutions are strongly affected by the number of used vehicles.
}

\begin{landscape}

\begin{table}[H]
\small
\centering
\caption{Summary of ROUTES-CONSTRUCTION executions and characteristics of the best solutions for instances from the first benchmark set.}
\begin{tabular}{ l | c c c | c c c c c c }
\hline
\multicolumn{1}{l|}{Instances} &  \multicolumn{3}{c|}{ROUTES-CONSTRUCTION executions} &  \multicolumn{6}{c}{characteristics of the best obtained solutions}\\
	&\#knpart	&\#bppart	&\#bpimpr	&min(veh)	&avg(veh)	&max(veh)	&min(col)	&avg(col)	&max(col)\\\hline
Fio-30acc-LOW-025	&4	&2	&0	&1	&1.0	&1	&4	&5.6	&9\\
Fio-30acc-MED-025	&64	&15	&0	&1	&1.3	&2	&3	&4.4	&8\\
Fio-30acc-HIGH-025	&83	&18	&0	&1	&1.4	&2	&4	&6.8	&11\\
Fio-60acc-LOW-025	&91	&57	&0	&2	&2.1	&3	&1	&3.7	&7\\
Fio-60acc-MED-025	&182	&107	&10	&2	&2.7	&3	&1	&3.9	&10\\
Fio-60acc-HIGH-025	&84	&58	&2	&2	&3.0	&4	&1	&2.8	&5\\
Dob-30acc-LOW-025	&2	&0	&0	&0	&0.6	&1	&5	&7.5	&9\\
Dob-30acc-MED-025	&0	&0	&0	&0	&0.9	&1	&6	&8.0	&10\\
Dob-30acc-HIGH-025	&0	&0	&0	&0	&0.9	&1	&5	&7.2	&12\\
Dob-60acc-LOW-025	&42	&23	&0	&1	&1.3	&2	&4	&5.9	&10\\
Dob-60acc-MED-025	&80	&28	&0	&1	&1.6	&2	&3	&5.7	&10\\
Dob-60acc-HIGH-025	&116	&52	&0	&1	&1.6	&2	&4	&5.4	&7\\
Fio-30acc-LOW-050	&50	&2	&0	&1	&1.1	&2	&2	&5.0	&8\\
Fio-30acc-MED-050	&63	&11	&0	&1	&1.4	&2	&2	&6.0	&10\\
Fio-30acc-HIGH-050	&61	&24	&0	&1	&1.4	&2	&1	&4.8	&9\\
Fio-60acc-LOW-050	&161	&122	&0	&2	&2.4	&4	&1	&5.2	&10\\
Fio-60acc-MED-050	&133	&68	&1	&2	&2.9	&3	&1	&4.3	&7\\
Fio-60acc-HIGH-050	&47	&26	&2	&2	&2.7	&4	&3	&5.1	&8\\
Dob-30acc-LOW-050	&8	&0	&0	&0	&0.7	&1	&2	&4.4	&6\\
Dob-30acc-MED-050	&0	&0	&0	&0	&0.9	&1	&3	&6.5	&12\\
Dob-30acc-HIGH-050	&0	&0	&0	&0	&0.9	&1	&6	&7.8	&10\\
Dob-60acc-LOW-050	&46	&30	&0	&1	&1.3	&2	&3	&4.6	&8\\
Dob-60acc-MED-050	&90	&8	&0	&1	&1.7	&2	&3	&6.5	&13\\
Dob-60acc-HIGH-050	&5	&1	&0	&1	&1.9	&3	&1	&3.4	&6\\
Fio-30acc-LOW-125	&47	&0	&0	&1	&1.1	&2	&2	&5.0	&8\\
Fio-30acc-MED-125	&75	&6	&0	&1	&1.4	&2	&3	&5.4	&8\\
Fio-30acc-HIGH-125	&32	&3	&0	&1	&1.6	&2	&3	&8.0	&13\\
Fio-60acc-LOW-125	&280	&117	&1	&2	&2.4	&3	&1	&4.0	&6\\
Fio-60acc-MED-125	&174	&129	&7	&2	&2.9	&4	&1	&4.4	&9\\
Fio-60acc-HIGH-125	&86	&55	&9	&2	&3.0	&4	&1	&2.5	&5\\
Dob-30acc-LOW-125	&4	&0	&0	&0	&0.7	&1	&2	&4.6	&9\\
Dob-30acc-MED-125	&0	&0	&0	&0	&0.9	&1	&3	&6.5	&9\\
Dob-30acc-HIGH-125	&0	&0	&0	&1	&1.0	&1	&3	&6.3	&11\\
Dob-60acc-LOW-125	&144	&1	&0	&1	&1.4	&2	&1	&4.6	&7\\
Dob-60acc-MED-125	&85	&13	&0	&1	&1.7	&2	&2	&3.5	&5\\
Dob-60acc-HIGH-125	&16	&1	&0	&1	&1.9	&2	&2	&3.9	&8\\ \hline
Averages	&65.4	&27.1	&0.9	&1.0	&1.6	&2.1	&2.6	&5.3	&8.7\\ \hline
	\end{tabular}
\label{tab:statsfio}
\end{table}

\begin{table}[H]
\small
\centering
\caption{Summary of ROUTES-CONSTRUCTION executions and characteristics of the best solutions for the small instances from the second benchmark set.}
\begin{tabular}{ l | c c c | c c c c c c }
\hline
\multicolumn{1}{l|}{Instances} &  \multicolumn{3}{c|}{ROUTES-CONSTRUCTION executions} &  \multicolumn{6}{c}{characteristics of the best obtained solutions}\\
	&\#knpart	&\#bppart	&\#bpimpr	&min(veh)	&avg(veh)	&max(veh)	&min(col)	&avg(col)	&max(col)\\\hline
20n-270r-2.5p	&0	&0	&0	&0	&0.6	&1	&2	&3.8	&5\\
20n-410r-2.5p	&6	&0	&0	&0	&0.9	&1	&2	&4.8	&7\\
20n-540r-2.5p	&0	&0	&0	&1	&1.0	&1	&2	&4.3	&7\\
20n-270r-3.5p	&1	&0	&0	&0	&0.6	&1	&2	&3.5	&5\\
20n-410r-3.5p	&7	&0	&0	&0	&0.9	&1	&3	&4.0	&5\\
20n-540r-3.5p	&0	&0	&0	&1	&1.0	&1	&2	&4.0	&6\\
25n-320r-2.5p	&1	&0	&0	&0	&0.7	&1	&2	&3.6	&7\\
25n-480r-2.5p	&5	&0	&0	&1	&1.0	&1	&4	&6.1	&9\\
25n-640r-2.5p	&34	&0	&0	&1	&1.3	&2	&2	&4.4	&6\\
25n-320r-3.5p	&0	&0	&0	&0	&0.7	&1	&3	&4.4	&5\\
25n-480r-3.5p	&7	&0	&0	&1	&1.0	&1	&2	&4.3	&8\\
25n-640r-3.5p	&46	&0	&0	&1	&1.3	&2	&2	&4.8	&8\\
30n-420r-2.5p	&6	&0	&0	&0	&0.9	&1	&1	&3.7	&5\\
30n-630r-2.5p	&11	&3	&0	&1	&1.1	&2	&3	&7.8	&13\\
30n-840r-2.5p	&44	&3	&0	&1	&1.6	&2	&2	&4.9	&10\\
30n-420r-3.5p	&0	&0	&0	&0	&0.9	&1	&3	&4.8	&7\\
30n-630r-3.5p	&14	&5	&0	&1	&1.1	&2	&4	&6.5	&10\\
30n-840r-3.5p	&84	&9	&0	&1	&1.6	&2	&3	&7.0	&12\\
35n-480r-2.5p	&3	&0	&0	&1	&1.0	&1	&4	&6.1	&11\\
35n-710r-2.5p	&16	&3	&0	&1	&1.3	&2	&4	&7.0	&10\\
35n-950r-2.5p	&116	&6	&0	&1	&1.9	&2	&1	&4.1	&8\\
35n-480r-3.5p	&6	&0	&0	&1	&1.0	&1	&5	&7.0	&9\\
35n-710r-3.5p	&10	&1	&0	&1	&1.3	&2	&5	&12.7	&16\\
35n-950r-3.5p	&88	&5	&0	&1	&1.9	&2	&2	&4.2	&7\\
40n-550r-2.5p	&0	&0	&0	&1	&1.0	&1	&6	&7.7	&10\\
40n-820r-2.5p	&164	&16	&0	&1	&1.6	&2	&1	&5.4	&8\\
40n-1090r-2.5p	&63	&0	&0	&2	&2.0	&2	&3	&5.1	&7\\
40n-550r-3.5p	&0	&0	&0	&1	&1.0	&1	&6	&9.6	&12\\
40n-820r-3.5p	&111	&12	&0	&1	&1.6	&2	&2	&5.6	&12\\
40n-1090r-3.5p	&56	&4	&0	&2	&2.0	&2	&4	&6.4	&9\\\hline
Averages &30.0	&2.2	&0.0	&0.8	&1.2	&1.5	&2.9	&5.6	&8.5\\ \hline
	\end{tabular}
\label{tab:statssmall}
\end{table}

\begin{table}[H]
\small
\centering
\caption{Summary of ROUTES-CONSTRUCTION executions and characteristics of the best solutions for the medium instances from the second benchmark set.}
\begin{tabular}{ l | c c c | c c c c c c }
\hline
\multicolumn{1}{l|}{Instances} &  \multicolumn{3}{c|}{ROUTES-CONSTRUCTION executions} &  \multicolumn{6}{c}{characteristics of the best obtained solutions}\\
	&\#knpart	&\#bppart	&\#bpimpr	&min(veh)	&avg(veh)	&max(veh)	&min(col)	&avg(col)	&max(col)\\\hline
50n-650r-2.5p	&65	&3	&0	&1	&1.3	&2	&2	&4.4	&7\\
50n-910r-2.5p	&49	&3	&0	&1	&1.7	&2	&5	&8.2	&14\\
50n-1300r-2.5p	&259	&37	&4	&2	&2.4	&3	&3	&5.6	&11\\
50n-650r-3.5p	&100	&2	&0	&1	&1.3	&2	&2	&5.0	&9\\
50n-910r-3.5p	&63	&4	&0	&1	&1.7	&2	&3	&6.3	&10\\
50n-1300r-3.5p	&186	&49	&2	&2	&2.4	&3	&1	&4.8	&7\\
60n-800r-2.5p	&148	&0	&0	&1	&1.6	&2	&3	&4.5	&8\\
60n-1190r-2.5p	&140	&11	&0	&2	&2.3	&3	&2	&3.3	&5\\
60n-1590r-2.5p	&21	&0	&0	&3	&3.0	&3	&2	&3.5	&8\\
60n-800r-3.5p	&96	&2	&0	&1	&1.6	&2	&1	&4.4	&7\\
60n-1190r-3.5p	&69	&3	&0	&2	&2.3	&3	&3	&5.8	&13\\
60n-1590r-3.5p	&42	&0	&0	&3	&3.0	&3	&2	&5.0	&9\\
80n-1070r-2.5p	&154	&2	&0	&2	&2.0	&2	&3	&4.7	&6\\
80n-1610r-2.5p	&147	&12	&3	&3	&3.0	&3	&2	&4.0	&8\\
80n-2150r-2.5p	&28	&3	&3	&4	&4.0	&4	&2	&4.6	&10\\
80n-1070r-3.5p	&217	&9	&0	&2	&2.0	&2	&3	&5.3	&8\\
80n-1610r-3.5p	&175	&18	&1	&3	&3.0	&3	&2	&4.3	&8\\
80n-2150r-3.5p	&35	&3	&3	&4	&4.0	&4	&2	&4.7	&8\\
100n-1330r-2.5p	&105	&21	&5	&2	&2.4	&3	&3	&6.1	&11\\
100n-2000r-2.5p	&181	&52	&18	&3	&3.7	&4	&2	&4.7	&10\\
100n-2670r-2.5p	&98	&3	&3	&4	&4.9	&5	&2	&5.0	&11\\
100n-1330r-3.5p	&133	&23	&6	&2	&2.4	&3	&4	&5.8	&10\\
100n-2000r-3.5p	&168	&48	&14	&3	&3.7	&4	&2	&4.5	&11\\
100n-2670r-3.5p	&77	&6	&6	&4	&4.9	&5	&2	&4.7	&9\\\hline
Averages 	&114.8	&13.1	&2.8	&2.3	&2.7	&3.0	&2.4	&5.0	&9.1\\ \hline
	\end{tabular}
\label{tab:statsmedium}
\end{table}

\begin{table}[H]
\small
\centering
\caption{Summary of ROUTES-CONSTRUCTION executions and characteristics of the best solutions for the large instances from the second benchmark set.}
\begin{tabular}{ l | c c c | c c c c c c }
\hline
\multicolumn{1}{l|}{Instances} &  \multicolumn{3}{c|}{ROUTES-CONSTRUCTION executions} &  \multicolumn{6}{c}{characteristics of the best obtained solutions}\\
	&\#knpart	&\#bppart	&\#bpimpr	&min(veh)	&avg(veh)	&max(veh)	&min(col)	&avg(col)	&max(col)\\\hline
120n-1410r-2.5p	&231	&48	&10	&2	&2.6	&3	&4	&7.8	&14\\
120n-2030r-2.5p	&154	&40	&14	&3	&3.7	&4	&2	&5.0	&9\\
120n-2820r-2.5p	&287	&158	&61	&5	&5.1	&6	&4	&8.7	&14\\
120n-1410r-3.5p	&252	&80	&12	&2	&2.6	&3	&4	&6.4	&9\\
120n-2030r-3.5p	&189	&35	&9	&3	&3.7	&4	&3	&6.0	&11\\
120n-2820r-3.5p	&245	&114	&60	&5	&5.1	&6	&3	&6.3	&11\\
160n-1810r-2.5p	&91	&10	&2	&3	&3.3	&4	&4	&8.3	&16\\
160n-2930r-2.5p	&119	&19	&13	&5	&5.4	&6	&2	&4.7	&7\\
160n-3320r-2.5p	&112	&59	&30	&6	&6.1	&7	&6	&12.4	&21\\
160n-1810r-3.5p	&140	&23	&7	&3	&3.4	&4	&2	&5.0	&8\\
160n-2930r-3.5p	&175	&50	&44	&5	&5.4	&6	&2	&3.8	&7\\
160n-3320r-3.5p	&105	&59	&25	&6	&6.1	&7	&1	&4.4	&10\\
200n-2110r-2.5p	&168	&33	&14	&3	&3.9	&4	&3	&5.4	&10\\
200n-3430r-2.5p	&189	&70	&44	&6	&6.3	&7	&3	&4.5	&9\\
200n-4120r-2.5p	&119	&47	&39	&7	&7.6	&8	&1	&3.6	&7\\
200n-2110r-3.5p	&119	&49	&19	&3	&3.9	&4	&3	&6.6	&11\\
200n-3430r-3.5p	&161	&59	&33	&6	&6.3	&7	&2	&4.3	&9\\
200n-4120r-3.5p	&105	&17	&11	&7	&7.6	&8	&2	&4.6	&9\\
300n-2510r-2.5p	&224	&111	&14	&4	&4.6	&5	&2	&5.8	&10\\
300n-4030r-2.5p	&182	&123	&46	&6	&7.6	&9	&1	&3.3	&7\\
300n-4520r-2.5p	&161	&103	&43	&8	&8.3	&9	&1	&4.0	&10\\
300n-2510r-3.5p	&259	&46	&7	&4	&4.6	&5	&4	&9.0	&15\\
300n-4030r-3.5p	&154	&97	&48	&7	&7.4	&8	&1	&3.9	&8\\
300n-4520r-3.5p	&189	&109	&46	&8	&8.3	&9	&1	&5.2	&11\\ \hline
Averages &172.1	&65.0	&27.1	&4.9	&5.4	&6.0	&2.5	&5.8	&10.5\\ \hline
	\end{tabular}
\label{tab:statslarge}
\end{table}

\end{landscape}

\section{Concluding remarks}
\label{sec:concludingremarks}

In this paper, we tackled the selective and periodic inventory routing problem (SPIRP) for waste vegetable oil collection by applying a fast and effective MIP-based heuristic. 
This relevant problem in the context of reverse logistics was an open area for research due to the fact that the best performing solution procedure available in the literature consumes a lot of time to provide solutions. 

Computational experiments have demonstrated that our newly proposed MIP-based heuristic is very fast and effective when providing near optimal solutions for the considered instances, achieving considerable low gaps within seconds. 
The approach has the benefit of providing an \textit{a posteriori} quality guarantee, as it contains the resolution of a valid relaxation of the problem in its mechanism\rafaelC{, thus providing both a feasible solution and a lower bound}.
Besides, the MIP-based heuristic was able to obtain high quality solutions, better than those found by a state-of-the-art solution procedure for several of the tested instances, while using just a fraction of the time.
The most significant performance was obtained for the larger instances, as our approach could improve the best known result for all of them. 

The SPIRP considered in this paper contains a single depot and a single type of product/oil. Thus, an immediate extension to be considered is to reformulate and solve this problem in a multi-depot environment. Also, it would be interesting to study the situation in which there are multiple types of product/oil with one or multiple depots.  These are some research opportunities that could be explored in the near future.

\vspace{0.8cm}

{
\noindent \small 
\textbf{Acknowledgments:} Work of Rafael A. Melo was supported by the State of Bahia Research Foundation (FAPESB) and the Brazilian National Council for Scientific and Technological Development (CNPq). \rafaelC{The authors are thankful to the anonymous reviewers for the insightful comments which helped improve the quality of this paper.}
}

\bibliographystyle{apacite}
\bibliography{main}
 
\end{document}